\documentclass[10pt,letterpaper]{article}
\usepackage[top=0.85in,left=2.75in,footskip=0.75in]{geometry}

\usepackage{amsmath,amssymb}

\usepackage{changepage}

\usepackage[utf8x]{inputenc}

\usepackage{textcomp,marvosym}

\usepackage{cite}

\usepackage{nameref,hyperref}

\usepackage{microtype}
\DisableLigatures[f]{encoding = *, family = * }

\usepackage[table]{xcolor}

\usepackage{array}

\newcolumntype{+}{!{\vrule width 2pt}}

\newlength\savedwidth

\raggedright
\usepackage[aboveskip=1pt,labelfont=bf,labelsep=period,justification=raggedright,singlelinecheck=off]{caption}

\makeatletter
\renewcommand{\@biblabel}[1]{\quad#1.}
\makeatother

\usepackage{multirow}

\usepackage{lastpage,fancyhdr,graphicx}
\usepackage{epstopdf}
\fancyheadoffset[L]{2.25in}
\fancyfootoffset[L]{2.25in}
\begin{document}
\vspace*{0.2in}

\begin{flushleft}
{\Large
\textbf\newline{Solving the Initial Value Problem of Ordinary Differential
Equations by Lie Group based Neural Network Method}  
}
\newline
\\
Ying Wen\textsuperscript{1\Yinyang},
Temuer Chaolu\textsuperscript{2*\Yinyang},
Xiangsheng Wang\textsuperscript{3},
\\
\bigskip
\textbf{1} College of Information Engineering, Shanghai Maritime
University, 201306, Shanghai, China
\\
\textbf{2} College of Arts and Sciences, Shanghai Maritime
University, 201306, Shanghai, China
\\
\textbf{3} College of Information Engineering, Shanghai Maritime
University, 201306, Shanghai, China
\\
\bigskip

%
%
\Yinyang These authors contributed equally to this work.





* tmchaolu@shmtu.edu.cn

\end{flushleft}
\section*{Abstract}
To combine a feedforward neural network (FNN) and Lie group (symmetry) theory of differential equations (DEs), an alternative artificial NN approach is proposed to solve the initial value problems (IVPs) of ordinary DEs (ODEs). Introducing the Lie group expressions of the solution, the trial solution of ODEs is split into two parts. The first part is a solution of other ODEs with initial values of original IVP. This is easily solved using the Lie group and known symbolic or numerical methods without any network parameters (weights and biases). The second part consists of an FNN with adjustable parameters. This is trained using the error back propagation method by minimizing an error (loss) function and updating the parameters. The method significantly reduces the number of the trainable parameters and can more quickly and accurately learn the real solution, compared to the existing similar methods. The numerical method is applied to several cases, including physical oscillation problems. The results have been graphically represented, and some conclusions have been made.



\section{Introduction}
Various fields, such as science, finance, and engineering, can transform several problems into a set of ordinary differential equations (ODEs) or partial DEs (PDEs) through mathematical modeling. Usually, the analytical solutions of these equations are unavailable. Therefore, solving the numerical solution of ODEs becomes particularly important to explore real world problems. Several off-the-shelf methods (e.g., Euler \cite{Rivertz2013}, Runge-Kutta \cite{Kaps1979}, Finite difference \cite{Noye2010}, Adomian decomposition \cite{Yun2019}, and Lie group \cite{Olver} methods) are available for solving numerical solutions of DEs \cite{Press}. Several new studies are being conducted to obtain more efficient algorithms. A class of them is the methods based on artificial neural networks (ANNs) or deep learning models \cite{Lagaris1998, Effati2010, ASH2000, Liying2007, Tsoulos2009, Beidokhti2009,Chakraverty2014, Mall2016, Dockhorn2019, LiS2020, eweinan, Tongyong, aa1, Dis}. The main idea used in these methods is the use of the highly accurate approximation capability of an ANN to a continuous function \cite{Tongyong}. This has been used in various ways. Considering \cite{Lagaris1998, Effati2010}, a representative ANN method for solving ODEs and PDEs is presented. Regarding the method, the trial solution is expressed as the sum of two terms. The first term satisfies the initial or boundary values, and it does not contain network parameters. The second term is a feedforward NN (FNN) to be trained to satisfy the given ODEs or PDEs. However, the first term is constructed synthetically by observing the given initial or boundary values. Such free selection of the first term in the trial solution may not be suitable for capturing nonlinearity information of the real solution in the domain of the initial point. Regarding \cite{ASH2000}, He et al. used the extended back propagation algorithm to train the derivative of an FNN to solve a class of first-order PDEs. Moreover, Li-ying et al. proposed an ANN algorithm based on cosine basis functions to solve ODEs \cite{Liying2007}. Ioannis et al. created a trial solution in the form of a neural network based on a syntactic evolution scheme to solve ODEs and PDEs \cite{Tsoulos2009}. Furthermore, Shekari et al. determined the approximate solutions of time-dependent PDEs based on ANNs and a hybrid method of minimization techniques and configuration methods \cite{Beidokhti2009}. Chakraver et al. proposed a regression-based NN model to solve the initial or boundary value problems (IBVPs) for ODEs \cite{Chakraverty2014}. Additionally, Mall et al. proposed a new method based on a single-layer Legendre NN model to solve the IBVP for nonlinear ODEs \cite{Mall2016}. Dockhorn proved that small NNs (parameters \textless 500) can accurately learn the solutions of practical problems (Poisson and Navier-Stokes equations) \cite{Dockhorn2019}. Shangjie obtained the approximate solutions of DEs through an ANN. In \cite{eweinan}, a deep learning algorithm is provided to solve the IBVP of a class of high dimensional stochastic PDEs. Multi-layer physical information NN deep learning was used to investigate data-driven peakon solutions and periodic peakon solutions of Camassa-Holm (CH) equation, Degasperis-Procesi equation, etc, with initial conditions \cite{2021Data}.

The goal of using ANN to solve ODEs or PDEs is to make it computationally cheaper. Generally, performing ANNs to solve ODEs or PDEs involves a number of parameters to be trained with a lot of data. Considering several DEs, there is prior knowledge that is currently not being used in ANN algorithms or machine learning practice. These include some mathematical principles, solution expressions, or physical information. This prior information can act as a regularization agent that reduces the number of trainable parameters or the demand for a large number of training data. Therefore, using such useful information in a network algorithm results in amplifying the information content of the solution that the algorithm guides itself towards the accuracy approximation, even when only small network models or a few training samples are available. The physics-inspired NNs for solving DEs are proposed, where physical conservation laws and prior physical knowledge are encoded into the NNs (refer to \cite{aa1} and references therein). For example, the logarithmic nonlinear Schr\"{o}dinger (LNLS) equation is solved by a deep learning method of physical information NN. This is an important physical model in several fields such as quantum optics and nuclear physics in \cite{2020Solving}. On the other hand, because the solution of DE is a continuous function, it is also important to find a trial function expression that is closer to the properties of the true solution. This makes the ANN effective and can quickly capture the implied properties of the real solution.

Lie group (symmetry) of DEs, universally known as one parameter transformation group method of DEs, has had a profound impact on all areas of mathematics (both pure and applied), physics, engineering, and other mathematically based sciences \cite{Olver, Bluman, Ibragimov}. The groups can be found using symbolic computational method, and they are used to construct the explicit solutions of the corresponding DEs \cite{chao}. They also reduce the dimension in order of equations or number of involved variables, making the method efficiently solve nonlinear DEs within a certain range. Particularly, regarding an IVP of a first-order system of ODEs, this method also provides an explicit formal formula of the solution. This point is used in the present study. The advantages of the Lie groups in numerical analysis are extensively discussed in \cite{Grobner, Shokin, Doron}.

This study presents a method for solving the IVP of ODEs by combining the Lie group theory of ODEs and an FNN, which is significantly computationally cost effective. Moreover, the method leads to a differentiable, closed analytic form solution of the problem on whole interval of independent variables. Considering the method, using Lie group theory, the trial solutions to the IVP are expressed by the sum of two parts. The first part is the solution to an IVP of other ODEs with the initial value of the original problem. The ODEs in the IVP are from the part terms of the original ODEs selected by obeying some solvable principles so that the first part of the solution can easily be solved in advance by Lie group method, some known symbolic or numerical methods. Owing to the specific feature of the IVP, the first part itself can capture the essential properties of the real solution in an interval of initial point of independent variable. Moreover, the first part contains no training parameters because it does not take part in training the networks. The second part of the solution is constructed by an FNN, without being required to satisfy any specific initial value.

The study makes the following contributions. This is first time that the Lie theory and ANN method are being combined to solve IVP of ODEs. Because the first part of the solution bears part of the workload and detects the nonlinearity of solution, small-scale networks and a small amount of training data (samples) can solve the IVP more accurately. The algorithm significantly increases efficiency of the network method to solve the IVP of ODEs. On the other hand, the method provides new idea and enlighten on how to add "complementary elements" and design an ANN to increase the efficiency of ANN algorithm to solve DEs by incorporating with mathematical theory. Furthermore, the method can be applied to more a general form IVP or BVP of ODEs. Moreover, after modifying or combining some of the existing methods \cite{Sun}, the approach can be applied to IVP and BVP of PDEs.

The rest of the paper is structured as following. In Section $2$, we introduce the Lie group method for solving IVP of ODEs and an expression of solution to the IVP is proved. In Section $3$, the general process of the proposed method is given. In Section $4$, numerical experiments on the application of our method to some problems are given, including the oscillation model in physical problems. Finally, in Section $5$, some discussions as well as future research directions are presented.

\section{Lie group expression for solution to the IVP of a DEs}
\subsection{A basic Theorem}
Let $G_{1}$ = $\left\{ T_{a} \right\}$ be a one parameter (denoted by $a$) transformations $ x^{*}=T_{a}(x)$ from $x =\left(x_{1}, x_{2}, \cdots, x_{n}\right) \in \mathbb{R}^{n}$ to $x^{*}=\left(x_{1}^{*}, x_{2}^{*}, \cdots, x_{n}^{*}\right) \in \mathbb{R}^{n}$ with components
\begin{equation}
T_{a}: x_{i}^{*}=f_{i}(x;a), i=1,2, \ldots, n. \label{transformation}
\end{equation}
Here functions $f =\left(f_{1}, f_{2}, \cdots, f_{n}\right)$ are continuous and differentiable with respect to (w.r.t) independent variables $x_i$ and parameter $a$ varying in a neighborhood $\mathcal{O}$ of $0\in \mathbb{R}^1$.

Let
\begin{equation}
D=\sum_{i=1}\xi_i(x)\partial x_i\label{operator1}
\end{equation}
be a linear partial differential operator generated from the auxiliary functions defined by
\begin{equation}
\xi^{i}(x)=\left.\frac{\partial f_{i}(x;a)}{\partial a}\right|_{a=0}, i=1,2, \ldots, n. \label{auxiliaryf}
\end{equation}

From \cite{Ibragimov}, we know that set $G_1$ form an one parameter transformation group with generator $D$ from $\mathbb{R}^n$ to itself if $f_i$ satisfy additional conditions: i) $f_i(x, 0)=x_i$; ii)  $f_i(f(x, a), b)=f_i(x, a+b)$; iii) $f_i$ is differential w.r.t parameter $a$. In the group case, by Taylor expansion at $a=0$, we have expressions of $x^*_i$ as follows
\begin{equation}
x_i^*=f_{i}(x;a)=x_{i}+\xi_{i}(x) a+O\left(a^{2}\right)=e^{a D} x_i,i=1,2,\cdots, n. \label{auxiliaryff}
\end{equation}
Here formal notation $e^{aD}=\sum_{k=0}^\infty\frac{a^k}{k!}D^k$ is used.

The correspondence between set $G_1$ and an one parameter transformation group is given in the following theorem (The more details can be seen in \cite{Olver, Bluman, Ibragimov}).

\textbf {Theorem  1} The set $G_1$ with functions $x_i^{*}=f_i(x ; a)$ in (\ref{transformation}) defining an one parameter transformation group satisfies the IVP of ODEs
\begin{equation}
\frac{d x^{*}}{d a}=D x^{*},\left.x^{*}\right|_{a=0}=x. \label{ivp}
\end{equation}
Conversely, for any given continuously differentiable functions $\xi_{i}(x)$, i.e., the operator $D$ determined by (\ref{operator1}), the function $x^{*}=f(x; a)$ obtained by solving the problem (\ref{ivp}) determines an one parameter transformation group with $D$ as generator. Moreover, the transformations $T_a: x^{*}=f(x; a)$, i.e. the solution of (\ref{ivp}) can be expressed as formula (\ref{auxiliaryff}), i.e., $x^*=e^{a D}x$.

\subsection{An expression of solution to IVP of an ODEs}
Effectively finding solution of IVP (\ref{ivp}) depends on the complexity of the operator $D$. To reduce this difficulty, let's assume that $D$ in (\ref{ivp}) can be decomposed into summands of two parts as
\begin{equation}
D=D_{1}+D_{2}. \label{operator2}
\end{equation}
Thus, we have

\textbf {Theorem  2}\ (A solution expression) For an operator $D$ in (\ref{ivp}) with decomposition part $D_1=\sum_{i=1}^ng_i(x)\partial_{x_i}$ in (\ref{operator2}) and initial values $x\in \mathbb{R}^n$, the solutions to (\ref{ivp}) have expression as
\begin{equation}
x^*=e^{a D} x=\bar{x}(x; a)+\tilde{\mathcal{N}}(x; a),  \label{Grobner}
\end{equation}
where $\bar{x}=\bar{x}(x; a)=e^{a D_{1}} x$ and $\tilde{\mathcal{N}}(x; a)=\int_{0}^{a} D_{2}\left(e^{(a-\tau) D} x\right)|_{x \rightarrow \bar{x}(\tau ; x)} d \tau $.

The proof of the theorem is obtained by Grobner's formula  in \cite{Grobner} given by
\begin{equation*}
e^{a D}=e^{a D_{1}}+\sum_{\alpha=1}^{\infty} \sum_{k=\alpha}^{\infty} \frac{a^{k}}{k !} D_{1}^{k-\alpha} D_{2} D^{\alpha-1}.
\end{equation*}

Obviously, the first part $\bar{x}=e^{a D_1}x$ of the above formula (\ref{Grobner}) is the solution to new IVP
\begin{equation}
\frac{d \bar{x}}{d a}=D_1 \bar{x},\left. \bar{x}\right|_{a=0}=x, \label{ivp1}
\end{equation}
by Theorem 1.

Evidently, by suitably choosing the first part $D_1$ of $D$ in (\ref{operator2}), IVP (\ref{ivp1}) is more easily solved than original IVP (\ref{ivp}) by Lie group method as well as various symbolic or numerical methods \cite{Press}. Moreover, in \cite{Ibragimov, Grobner}, it is proved that under probably selection $D_1$, the solution of (\ref{ivp1}) higher accurately approximates the real solution of (\ref{ivp}) in some interval of initial point $a=0$.

Although the second part $\tilde{\mathcal{N}}$ in (\ref{Grobner}) will computed by symbolically in some simple cases, in general it is extremely complicate if one directly solve it. In this article, we use the formula (\ref{Grobner}) to design a FNN to solve the term.

\textbf{Remark 1:} It is notice that from (\ref{Grobner}), we have  $\tilde{\mathcal{N}}(x, 0)=0$. Hence without loss of generality we suppose $\tilde{\mathcal{N}}(x; a)=a\bar{\mathcal{N}}(x; a)$ for a differential function $\bar{\mathcal{N}}(x; a)$.

\section{Method}

In the section, we describe the scheme of our Lie group based on FNN method which relies on the expressions (\ref{Grobner}) of the solutions to an IVP of ODEs given in Theorem 2.
\subsection{Scheme}
Our method, first, is designed for the IVP of an autonomous system of ODEs
\begin{eqnarray}
&&\frac{dy_i}{dx}=f_i(y_1,y_2,\cdots, y_n),\,\,y_i(0)=\alpha_i\in \mathbb{R}^1, i=1, 2, \cdots, n,\label{equation}
\end{eqnarray}
where $x\in \mathcal{O}\subset \mathbb{R}^1$ is independent variable and $y_i=y_i(x)$ are dependent variables and $f_i$ are differential functions of own arguments. Then, we extend the method to solve different form ODEs or PDEs problems.

We rewrite IVP (\ref{equation}) as operator form
\begin{equation}
\frac{dy}{dx}=D y, \,\, y(0)=\alpha,\label{eq1}
\end{equation}
by introducing notations $y=(y_1, y_2, \cdots, y_n)$, $\alpha=(\alpha_1, \alpha_2,$ $ \cdots, \alpha_n)$ and operator $D=\sum_{i=1}^n f_i(y)\frac{\partial}{\partial y_i}$, which has the same form with (\ref{ivp}). Consequently, by Theorem 2, the solution of IVP (\ref{equation})  can be written as
\begin{equation}
y=e^{x D} \alpha=e^{x D_{1}} \alpha+x \bar{\mathcal{N}}(x; \alpha)  \label{Grobner1}
\end{equation}
for a decomposition $D=D_1+D_2$ and functions $\bar{\mathcal{N}}(x; \alpha)=\{\bar{\mathcal{N}}_i(x, \alpha)\}_{i=1}^n$ with $x$ as group parameter. Here we understand that $e^{xD}\alpha=e^{xD}y|_{y\rightarrow\alpha}$.

We use this expression of the solution $y$ to design a FNN method to approximate $y(x)$ by a neural network $\hat{y}(x, \theta)$, where $\theta$ is the neural network parameters. To this end, we require that the network (trial solution) has the expression
\begin{eqnarray}
\hat{y}(x)=\bar{y}(x)+x\mathcal{N}(x;\theta),\label{trial}
\end{eqnarray}
where $\mathcal{N}(x;\theta)$ is a differential function to be determined and  $\bar{y}(x)=e^{xD_1}\alpha$ which is solution of the following new IVP
\begin{eqnarray}\frac{d\bar{y}}{dx}=D_1 \bar{y}=g(\bar{y}), \,\, \bar{y}(0)=\alpha, \label{part1}
\end{eqnarray}
by Theorem 1 for $D_1=\sum_{i=1}^ng_i(\bar{y})\partial_{\bar{y}_i}$.

Therefore, the deviations between the trial solution $\hat{y}$ in (\ref{trial}) and exact solution $y$ in (\ref{Grobner1}) is characterized by
\begin{eqnarray}\Delta y=|x(y(x)-\hat{y}(x))| \label{deviation}\end{eqnarray}
for $x$ in considered interval.

\textbf{Remark 2: }We call IVP (\ref{part1}) as an associated IVP of original IVP (\ref{equation}) or (\ref{eq1}). Its solution $\bar{y}(x)$ is determined in advance by solving IVP (\ref{part1}) and by applying the Lie group method or other acknowledged various methods after appropriately selecting operator $D_1$. The solvability and good approximation capability of the associated IVP solution $\bar{y}$ are main considerations of choosing operator $D_1$.

Consequently, compared with synthetic ways to give the first part of the trial solution in literatures, here the first part of the trial solution is constructed by a solution of the associated IVP which more accurately approximates the real solution of IVP (\ref{equation}) to be solved in an interval of initial point $x=0$. Due to the promising feature of efficient approximation of the first term in (\ref{trial}) to exact solution, the construction of the second term in (\ref{trial}) by networks algorithm takes much less a computational burden than that of solving the overall solution with big scale training parameters. Hence, in automatic learning (training) procedure, the network detecting quickly tunes itself to approach the natural properties of the solution. This is the main our motivation why we combine the Lie group method with ANN. Therefore, our attention focuses on determining the second part of solution (\ref{trial}).

This unknown function $\mathcal{N}(x; \theta)$ in the second part of the trial solution (\ref{trial}) is constructed by a FNN through optimizing the loss (error) function
\begin{equation}
\mathcal{L}(\theta)=\frac{1}{N}\sum_k^N\sum_{i=1}^n\{\frac{d\hat{y}^k_i}{dx}-
f_i(\hat{y}^k_1,
\hat{y}^k_2,\cdots, \hat{y}^k_n)\}^2,\label{loss}
\end{equation}
where $\hat{y}^k_i=\hat{y}_i(x_k, \theta),\frac{d\hat{y}^k_i}{dx}=\frac{d\hat{y}^k_i (x,\theta)}{dx}|_{x\rightarrow x_k}$ and the data set $S=\{x_k\}_{k=1}^N$ is a set of train points obtained from training interval $\mathcal{O}$ in some distributive sense.

Assuming that the operator $D=D_1+D_2$ in (\ref{eq1}) has expressions
\begin{eqnarray}
&&D_1=\sum_{i=1}^ng_i(y)\partial y_i,\,\,D_2=\sum_{i=1}^nh_i(y)\partial y_i \textrm{ with } f_i=g_i+h_i, \label{D12}
\end{eqnarray}
Thus the loss function (\ref{loss}) is rewritten as
\begin{equation}
\mathcal{L}(\theta)=\frac{1}{N}\sum_k^N\sum_{i=1}^n\left(g_i(\bar{y}^k)+(\frac{d (x N_i(x;\theta))}{dx}|_{x=x_k}-
f_i(\hat{y}^k))\right)^2,\label{loss1}
\end{equation}
obtained by substituting (\ref{trial}),(\ref{part1}) and (\ref{D12}) into (\ref{loss}).

In the loss function (\ref{loss1}), the $g_i(\bar{y}^k)$ is independent of training parameters $\theta$. This saves more computational effort. This is another benefit of combining Lie expression of DE solution and ANN algorithm.

In our case (\ref{equation}), the network models should have one- input $x$ and an $n$- output $\mathcal{N}(x;\theta)=\{\mathcal{N}_i(x;\theta)\}_{i=1}^n$ layer. We use the network  possessing $m$ neurons as model to present our proposed method scheme.

For each input $x$, there are outputs $\mathcal{N}_i(x;\theta)$ as follows
\begin{eqnarray}
&&\mathcal{N}(x;\theta)=\delta(\hat{\mathcal{N}}_1,\hat{\mathcal{N}}_2,\cdots, \hat{\mathcal{N}}_n,),\,\,\, \hat{\mathcal{N}}_j=\sum_{i=1}^{m} v_{ji} \cdot \sigma(x \cdot w_{i}+b_{i})+c_j, \label{N}
\end{eqnarray}
where trainable parameters $\theta=\{w_{i}, v_{ji}, b_k, c_j\}$ and $w_i$ is the weight from the input $x$ to the $i$th neuron in the hidden layer and $v_j=\{v_{j1}, v_{j2},\cdots, v_{jm}\}$ and $v_{ji}$ is the weight vector from the $i$th neuron in the hidden layer to the $j$th neuron in output layer, $b_{j}$ is the bias of the $j$th neuron in hidden layer and $c_j$ is the bias of the $j$th output neuron. The both $\sigma(\cdot)$ and $\delta(\cdot)$ are activation functions for outlets of hidden and output layers respectively. In our examples given next section, we use a Sigmoid function as hidden layer activation and take a linear activation in the output layer for our networks.

In training the networks, the values of loss function are obtained by the forward training with formula (\ref{N}), and the gradient descents of loss function (\ref{loss1}) are obtained in the backpropagation; minimizing the loss function involves not only the network output $\mathcal{N}(x; \theta)$, but also the derivatives of the network output $\mathcal{N}(x;\theta)$ w.r.t the input and network parameters; while these derivatives are given recursively by the ones of $\hat{\mathcal{N}}_i$ and functions $\delta$ and $\sigma$; the update network parameters are using unified formula as
\begin{eqnarray*}
\theta^{k+1}=\theta^{k}+\Delta \theta^{k}=\theta^{k}-\eta \frac{\partial \mathcal{L}(\theta)}{\partial \theta^{k}},
\end{eqnarray*}
where $\eta$ is the learning rate and $k$ is the iteration step (refer to \cite{Lagaris1998}).
\subsection{Principle to select operator $D_1$}

The efficiency of the proposed method above depends on the selecting of operator $D_1$ in decomposition $D=D_1+D_2$. The main principles in choosing $D_1$ are given in the following considerations.

1. The associated IVP (\ref{part1}) is easily solved explicitly or numerically. In particular, we choose $g_i$ in (\ref{D12}) as part terms of $f_i$ in (\ref{equation}) so that the associated IVP (\ref{part1}) yields explicit  or numerical solutions.

2. The operations of calculation on loss function (\ref{loss1}), such as calculation of derivatives and values of them, are not expansive as possible as.
\subsection{Extension}
The method proposed above can be applied to any IVP or BVP of DEs which as long as can be transformed to the form of (\ref{equation}).

\subsubsection{Non-autonomous case.}
For the non-autonomous case
\begin{eqnarray}
&&\frac{dy_i}{dx}=f_i(x, y_1, y_2, \cdots, y_n),\,\, y_i(0)=\alpha_i\in \mathbb{R}^1, i=1, 2, \cdots, n,\label{nonautoequation}
\end{eqnarray}
we introduce a new variable $y_0=x$ and add an equation $\frac{dy_0}{dx}=1$ and initial value $y_0(0)=0$ to (\ref{nonautoequation}). Then, the current IVP turns to the standard IVP (\ref{equation}) with $n+1$ unknown functions $y_0, y_1, \cdots, y_n$ and operator $D=\sum_{i=1}^nf_i\partial_{y_i}+\partial_x$.
\subsubsection{High-order case.} For the general  IVP of an $n$-order non-autonomous ODE
\begin{eqnarray}
&&\frac{d^{n} y(x)}{d x^{n}}=f(x, y, y', \ldots, y^{(n-1)}),\nonumber \\
&&y^{(k)}=a_{k+1}, \,\,\, k =0, 1, \ldots, n-1, \label{korder}
\end{eqnarray}
we can transform it to the standard IVP (\ref{equation}) by introducing transformations $y_0=x, y_1=y, y_2=y_1^\prime, \ldots, y_{n}=y_{n-1}^\prime$ with correspondingly $f_0=1, f_1=y_2, \ldots, f_{n-1}=y_n, f_n=f(y_0, y_1,$ $\ldots,y_{n})$.
\subsubsection{nonzero start points}
In nonzero start point $x=x_0$ in (\ref{ivp}), we do a translation $x\rightarrow x-x_0$ to change the problem into standard form.
\subsubsection{PDE case}
The most intuitive idea is to convert the PDE to the ODE and then solve the original problem indirectly. For example, use symmetries of PDE, traveling wave transformations, etc to decrease the order of the PDEs and number of independent variables involved in the PDEs. Another idea is to combine existing semi-discrete methods for solving PDE \cite{Sun}. Discrete the underline PDEs in spatial variables approximately so that the resulting semi-discrete equation can be cast into an ODEs (\ref{equation}), then use the proposed method on the resulted ODEs.
\subsubsection{Two point problem case}
To two points problem of ODEs, we firstly regard the problem as an IVP with using a point value as initial value problem, then put another point value in the loss function so that the network automatically learn the boundary value.

Next section, by applying the given method on solving some physical problems, we illustrate the superiority of the method.

\section{Numerical Experiments and Applications}
In this section, we provide some examples to show the effectiveness of our proposed method. Considering each example, we evaluate the method by calculating the error between our results and the exact solutions (if available) or numerical solutions obtained from known numerical methods. Regarding all the examples, a small architecture FNN with one hidden layer of only three neurons, a linear output, and a small training data is used. To illustrate the feature of the solution obtained by our method, we provide figures displaying the graph of the solution and the exact solution of the training interval. In addition, we draw the compared figures of first terms of the trial solutions given in present article and literatures to show the capabilities of the terms to capture nonlinear properties of the real solutions. Moreover, we consider points outside the training interval (some extensions of the training interval) to show the generalization and stability of our method.

\textbf{Example 1. }The first example shows the procedure and efficiency of our proposed method by comparing it with the methods provided in previous studies. We solve the IVP of two coupled first-order nonlinear ODEs considered in \cite{Lagaris1998}.

\begin{eqnarray}
\begin{aligned}
&\frac{d y_{1}}{d x}=\cos (x)+y_{1}^{2}+y_{2}-\left(1+x^{2}+\sin ^{2}(x)\right), \\
&\frac{d y_{2}}{d x}=2 x-\left(1+x^{2}\right) \sin (x)+y_{1} y_{2},
\end{aligned} \label{equation10}
\end{eqnarray}
with $x\in [-1,1]$ (training interval) and initial values $y_{1}(0)=0$ and $y_{2}(0)=1$. The IVP has exact solutions $y_{1}(x)=\sin (x)$ and $y_{2}(x)=1+x^{2}$.

To apply our method, let's turn the problem into the form of (\ref{equation}), i.e., an autonomous system of equations
\begin{eqnarray}
\left\{
\begin{aligned}
&y'_{1}(x)=f_1(y_0, y_1, y_2)=\cos (y_{0})+y_{1}^{2}+y_{2}-\left(1+y_{0}^{2}+\sin ^{2}(y_{0})\right), \\
&y'_{2}(x)=f_2(y_0, y_1, y_2)=2 y_{0}-\left(1+y_{0}^{2}\right) \sin (y_{0})+y_{1} y_{2},\\
&y'_{0}(x)=f_0(y_0, y_1, y_2)=1,
\end{aligned}\right. \label{equationaut10}
\end{eqnarray}
with introducing variables $y_0=x$ and using initial values $y_0(0)=y_1(0)=0$, $y_2(0)=1$.

From (\ref{D12}), the differential operator $D$ of this equation is
$
D=f_1 \partial y_{1}+f_2\partial y_{2} +f_0\partial y_{0}$. We select the linear parts of functions $f_i$ to construct operator $D_{1}=g_1\partial y_{1}+g_2\partial y_{2}+g_0\partial y_{0}$ in which $g_1(y_0,y_1,y_2)=\cos({y}_0)+{y}_2-(1+{y}_0^2+\sin^2({y}_0))$, $
g_2=2{y}_0-(1+{y}_0^2)\sin({y}_0)$, $g_3=1$. Hence, the associated IVP of (\ref{equationaut10}) corresponding to operator $D_{1}$ is a linear system
\begin{eqnarray}
\frac{d\bar{y}_i}{dx}=e^{xD_1}\bar{y}_i=g_i(\bar{y}_0,\bar{y}_1,\bar{y}_2), i=0,1,2,  \label{equation10sys}
\end{eqnarray}
with initial values $\bar{y}_0(0)=\bar{y}_1(0)=0, \bar{y}_2(0)=1$. Its exact solutions $\bar{y}_0=x$ and
\begin{eqnarray*}
&&\bar{y}_{1}(x)=x^2 \sin (x)+x/2-4\sin (x)+\sin (2 x)/4+4x \cos (x), \\ &&\bar{y}_{2}(x)=2 +x^2+x^2 \cos (x)-2 x \sin (x)-\cos (x),
\end{eqnarray*}
are easily obtained by Computer algebra system Mathematica. Now, the trial solutions of the two neural networks are supposed to be
\begin{eqnarray}&&{\bar{y}_{1}}(x)=\bar{y}_1+x\mathcal{N}_{1}\left(x;\theta\right),\,\,\,
{\bar{y}_{2}}(x)=\bar{y}_2+x  \mathcal{N}_{2}\left(x;\theta\right).\label{exam0}
\end{eqnarray}

The first parts of the trial solutions are simply adopted as $\bar{y}_1=0$ and $\bar{y}_2=1$ \cite{Lagaris1998}. Considering Fig. \ref{fig1}, the degrees to which the various first parts $\bar{y}_{1,2}$ used in both our trial solutions (\ref{exam0}) and reference \cite{Lagaris1998} approximate the true solutions near the initial point $x=0$ are shown. The first parts used in (\ref{exam0}) efficiently capture the nonlinear properties of the real solutions. This accelerates the convergence of the subsequent neural computations in our method.

\begin{figure}[!h]
\caption{{\bf Comparison of the efficiencies of the first terms $\bar{y}_{1,2}$ in the trial solutions approximating the exact ones $y_{1,2}$ in an interval of $x=0$ in Example 1.}
A:Comparison of $\bar{y}_{1}$ and $y_{1}$ at the initial point $t=0$. B:Comparison of $\bar{y}_{1}$ and $y_{1}$ at the initial point $t=0$.}
\label{fig1}
\end{figure}


Our networks consist of one hidden layer with three neurons trained on a set of only $21$ equidistant points in the interval $[-1,1]$. They yield higher accuracy approximate solutions of IVP (\ref{exam0}) as shown in the left figure of Fig. \ref{fig2}. Moreover, the solutions obtained by our method can approximate the exact solutions very well outside the training set till at least interval $[-1.5, 1.5]$ as shown in the right figure of Fig. \ref{fig2}. This indicates that this method has a higher accuracy on the whole training interval, better generalization (inertia), and stability, even when using small-scale networks and data set for the training.

\begin{figure}[!h]
\caption{{\bf Comparison between our solutions $\hat{y}_{1,2}$ and the exact solutions $y_{1,2}$ in Example 1.}}
\label{fig2}
\end{figure}


Considering Fig. \ref{fig3}, the accuracies by computing the deviations $\Delta_i={y}_i-\hat{y}_i$ between the trial solutions $\hat{y}_i$ obtained by our method and the exact solutions $y_i$ on the training set are controlled in $10^{-5}$ order of magnitudes. This shows that our method admits strong robustness on the training interval by fewer training samples. Regarding the predicted (no training) intervals $[-1.5,-1]$ and $[1,1.5]$, although the errors are increasing, the accuracies are maintained in $10^{-3}$ order of magnitudes without the training data.

\begin{figure}[!h]
\caption{{\bf Accuracies of the $\hat{y}_{1, 2}$ approximating exact ones $y_{1,2}$ in Example 1.}}
\label{fig3}
\end{figure}


Regarding Fig. \ref{fig4}, the dependence of the network errors on the number of iteration steps in our methods is shown. Considering the figure, when the iteration is approximately $130$ step, the error is close to $0$, and this state is persisted. This shows that the convergence of our method is faster, demonstrating the stability of our algorithm.

\begin{figure}[!h]
\caption{{\bf Trend of training errors with iteration steps in Example 1.}}
\label{fig4}
\end{figure}


Compared to the method in a previous study \cite{Lagaris1998} in which network architecture consisted of one hidden with ten neurons, our method achieves higher accuracy solutions in a shorter computational time by using fewer network parameters and small scale training data.

\textbf{Example 2}. We consider the linearly forced oscillation problem
\begin{eqnarray}
\ddot{y}+2\epsilon \dot{y}+y=f(t),
y(0)=a,\dot{y}(0)=b, \label{ex1}
\end{eqnarray}
where $f(t)$ is a known function and $a$ and $b$ are constants and $\epsilon$ is damping parameter with $0<\epsilon<1$. Since the solution of the problem is sensitive to the parameter $\epsilon\rightarrow 0$. It is well known that it was solved by Multiple scales and averaging methods \cite{perturbation}. To apply our method, let's turn the problem into the form of (\ref{equation})
\begin{eqnarray}
&&\dot{y}_0=1,\,\,\, \dot{y}_1=y_2,\dot{y}_2+2\epsilon y_{2}+y_{1}=f(y_0), \label{exam1}
\end{eqnarray}
with initial values $y_0(0)=0, y_1(0)=a, y_2(0)=b$ by letting $y_0=t, y_1=y$, which yields the operator $D=y_2\partial_{y_1}+(f(y_0)-y_1-2\epsilon y_2)\partial_{y_2}+\partial_{y_0}$. In order to its associated IVP more accurately approximate the real solution of (\ref{exam1}) and be solved easily, we take $D_1=y_2\partial_{y_1}-(y_1+2\epsilon y_2)\partial_{y_2}+\partial_{y_0}$
and $D_2=f(y_0)\partial_{y_2}$ as in (\ref{D12}).

By symbolic computation, one obtains the exact solutions of the associated IVP of (\ref{exam1}) for selected the $D_1$ as
$\bar{y}=(\bar{y}_0, \bar{y}_1, \bar{y}_2)=e^{tD_1}(y_0, y_1, y_2)|_{y_0\rightarrow 0, y_1\rightarrow a, y_2\rightarrow b}$, i.e., $\bar{y}_0=t, \bar{y}_1=\int^t_0\bar{y}_2(t)dt$ with
\begin{equation*}\bar{y}_1={\sigma}^{-1}{e^{-\epsilon t}}\left((a \epsilon +b) \sin \left(\sigma t \right)+a \sigma \cos \left(\sigma t \right)\right), \,\,\, \textrm{with }\,\, \sigma=\sqrt{1-\epsilon^2}.
\end{equation*}
Therefore, we have trial solution as $$\hat{y}_1=\bar{y}_1+t \mathcal{N}(t,\theta).$$
To concretely computation, we take example $\epsilon=0.2, a=0, b=1$ and $f(t)=-0.4 e^{-0.4t}\cos t$. In this situation, IVP (\ref{ex1}) has exact solution
$$y(t)=e^{-0.4t} \sin \left(t\right)$$

Regarding the previous studies, the first part of the trial solutions should be simply considered $\bar{y}_1=0$ or $x$ or $x(1+x)$. The degrees to which these first terms of the trial solutions approximate the exact solution in the neighborhood of initial point are shown in the Fig. \ref{fig5}. This shows that the first term in the trial solution $\hat{y}_1$ more accurately captures the nonlinearity of the exact solution $y$ in an interval of the initial point.

\begin{figure}[!h]
\caption{{\bf Comparison of the efficiencies of first term $\bar{y}_1$ in trial solutions $\hat{y}_1$ approximating the exact one $y$ in an interval of $t=0$ in Example 2.}}
\label{fig5}
\end{figure}


We use a grid of $21$ equidistant points in $[0, 2]$ as the training and testing data for training our networks $\mathcal{N}(t, \theta)$. Using such a small-scale network and fewer training data, our method obtains more accurate results.

The comparisons between the exact solutions ${y}_{1,2}$ and approximated ones $\hat{y}_{1, 2}$ provided by our method on the training interval $[0, 2]$ and its extension $[0, 2.5]$ are shown in Fig. \ref{fig6}.

\begin{figure}[!h]
\caption{{\bf Comparison of the exact solutions $y_i$ and our $\hat{y}_i$ in Example 2.}}
\label{fig6}
\end{figure}


Considering Fig. \ref{fig7}, the trend of the network errors $\mathcal{L}(\theta)$ as the number of iteration steps increases is shown. The graph reflects the rapid decline in the algorithm, and the error is close to zero at approximately $220$ iterations, reaching $10^{-7}$. This indicates that the algorithm converges quickly and is stable.

\begin{figure}[!h]
\caption{{\bf Trend of training the errors with iteration steps in Example 2.}}
\label{fig7}
\end{figure}


Considering Fig. \ref{fig8}, the deviations $\Delta y_{i}$ of solution $y_i$ on the training and test points are shown. It can be observed that the algorithm has a good generalization performance and a higher accuracy even on few training samples and a small network structure.

\begin{figure}[!h]
\caption{{\bf Accuracies of the $\hat{y}_{1, 2}$ to exact ones $y_{1, 2}$ in Example 2.}}
\label{fig8}
\end{figure}

To show the robustness of our method, we investigate the performances of our method on different values of parameter $\epsilon$ in (\ref{ex1}) for the same nonhomogeneous term $f(t)$. The results are shown in the following Table \ref{table1}.

\begin{table*}[!htbp]
\begin{adjustwidth}{-2.25in}{0in} 
\centering
\caption{
{\bf Robustness of our method.}}
\begin{tabular}{|c|c|c|c|}
\hline
&\text{$\epsilon$}&\text {Ave.Error }&\text {Ave.Deviation } \\\hline
\multirow{4}*{\shortstack{Training Data\\(21 samples)}}&0.01 & $7.604\times 10^{-6}$& $1.958\times10^{-3}$\\ 
\cline{2-4}
&0.2& $4.736\times10^{-7}$& $9.122\times10^{-5}$ \\
\cline{2-4}
&0.5& $2.174\times10^{-6}$& $3.166\times10^{-5}$  \\
\cline{2-4}
&0.8 & $1.646\times10^{-6}$& $8.282\times10^{-5}$\\
\cline{1-4}
\multirow{4}*{\shortstack{Predict Data\\(21 samples)}}&0.01 & $6.671\times10^{-6}$& $9.546\times10^{-4}$\\
\cline{2-4}
&0.2& $2.507\times10^{-7}$& $6.366\times10^{-5}$ \\
\cline{2-4}
&0.5& $1.268\times10^{-6}$& $3.565\times10^{-5}$  \\
\cline{2-4}
&0.8 & $2.168\times10^{-6}$& $6.601\times10^{-5}$\\
\cline{1-4}
\multirow{4}*{\shortstack{Test Data\\(26 samples)}}&0.01 & $2.318\times10^{-5}$& $1.925\times10^{-3}$\\
\cline{2-4}
&0.2& $6.549\times10^{-5}$& $5.018\times10^{-4}$ \\
\cline{2-4}
&0.5& $2.806\times10^{-4}$& $1.409\times10^{-4}$  \\
\cline{2-4}
&0.8 & $1.922\times10^{-5}$& $1.870\times10^{-4}$\\
\hline
\end{tabular}
\label{table1}
\end{adjustwidth}
\end{table*}

From the table, we can see that the strong stability of our algorithm is presented.

\textbf {Example 3. }Consider the following nonlinear initial value problem of Duffing equation \cite{perturbation}
\begin{equation}
u''+u+2\epsilon u^3=0  \label{Duffing}
\end{equation}
with initial conditions $u(0)=1, u'(0)=0$ and $0<\epsilon<1$. In mechanics, it is solved by perturbation techniques, such as strightforward expansion, multiple scales, averaging methods, etc to analysis the different resonance phenomenons.

Here, as an application of our proposed method, we consider the case with parameter values $\epsilon=0.5$.

One of the standard forms (\ref{equation}) for equation (\ref{Duffing}) is
\begin{eqnarray}
\begin{aligned}
&y'_{1}(t)= y_{2}(t),\\
&y'_{2}(t)=-y_{1}(t)-y^3_{1}(t),
\end{aligned} \label{equation1au}
\end{eqnarray}
with initial values $y_1(0)=1, y_2(0)=0$ by introducing variables $y_1=u, y_2=\dot{u}$. Obviously, we have operator $D=y_{2}\partial y_{1}+ (-y_{1}(x)-y^3_{1}(x))\partial y_{2}$. As before, one splits the operator into the sum of two parts and select first part as $D_{1}=y_{2}\partial y_{1} -y_{1}\partial y_{2}$. It produces the associated IVP
\begin{eqnarray}
\dot{\bar{y}}_1=\bar{y}_2,\,\,\, \dot{\bar{y}}_2=-\bar{y}_1-\bar{y}_1^3, \,\,\, \bar{y}_1(0)=1, \bar{y}_2(0)=0,\label{Duff}
\end{eqnarray} and the trial solutions $\hat{y}_{1}(t)=\bar{y}_1+t \mathcal{N}_{1}\left(t; \theta\right)$ and $\hat{y}_{2}(t)=\bar{y}_2+t \mathcal{N}_{2}\left(t; \theta\right)$ of (\ref{equation1au}) where $\bar{y}_1$ and $\bar{y}_2$ are solutions of (\ref{Duff}). The IVP (\ref{Duff}) has an obvious exact solution. However, it corresponds to (\ref{Duffing}) which is easily solved by the known numerical methods, such as Runge-Kutta method. This numerical solutions approximate the real solution more well than that given in literatures at initial point $t=0$ shown in Fig. \ref{fig9}.

\begin{figure}[!h]
\caption{{\bf Comparison of the approximate efficiencies of first terms  $\bar{y}_1, \bar{y}_2$ to the exact solution $y_1, y_2$ at initial point t=0 in Example 3.}
A:Comparison of $\bar{y}_{1}$ and $y_{1}$ at the initial point $t=0$. B:Comparison of $\bar{y}_{2}$ and $y_{2}$ at the initial point $t=0$}
\label{fig9}
\end{figure}


By using the exact solutions as first part of the trial solutions, we train NNs for $\mathcal{N}_1(t;\theta)$ and $\mathcal{N}_2(t;\theta)$ with a grid of $21$ equidistant training points in $[0,2]$. Fig. \ref{fig10} compares the training results with the Runge-Kutta method. It shows that higher accuracy approximation and quicker convergence are made throughout the whole domain.

\begin{figure}[!h]
\caption{{\bf Comparison between our solutions $\hat{y}_i$ and Runge-Kutta results $y_i$ in Example 3.}}
\label{fig10}
\end{figure}


Considering Fig. \ref{fig11}, the decreasing trend of train errors as iteration steps increasing is displayed showing the quick convergence and stability of our algorithm.

\begin{figure}[!h]
\caption{{\bf Trend of training errors with iteration steps in Example 3.}}
\label{fig11}
\end{figure}


It shows that the proposed method can also capture the fluctuation wave properties of the solution to a strong nonlinear dynamics system.
\section{Discussion and Conclusions}
In this study, we propose a new NN method based on the Lie group theory of ODEs to solve the IVP of ODEs. Examples are used to verify the effectiveness of the method, and the superiority of the method is proven by comparing the exact solutions and other numerical methods. The examples also prove that even when the network structure is small, this method can achieve higher accuracy solutions than that of the methods in \cite{Lagaris1998}. Moreover, the combination of the Lie group and NN methods is easy to implement.

1. The success of the method can be attributed to two aspects. The first is the employment of the Lie group expression of the real solution to the IVP of ODEs. These expressions provide a two-part summation form of the solution. Considering the first part, the initial values can be easily determined before the network calculation stage and is independent of the trainable parameters. Moreover, the first part yields approximate solutions to the real solution to be determined at least in an interval of the initial point of the independent variable.
These features of the first part of the trial solution save much workload in machine learning, which results in higher efficiency of the method. The second is the use of NNs that are excellent function approximators to determine the second part of the solution. The training (learning) procedure is implemented by optimizing the loss functions derived from the original ODEs and the additional initial or boundary values satisfied by the trial solutions.

2. Contrary to most previous methods, the proposed method is more general and can be extended to solve various problems of ODEs and PDEs by the appropriate selection of Lie group expressions of trial solutions and the loss functions.

3. As indicated by the applications in the numerical experiments, the method exhibits excellent generalization and stability performance. This is because the deviations in the training and testing data are of lower values and are uniformly stable.

4. The NN architectures employed were one hidden layer with three neurons, indicating that the trainable parameters are fewer. This indicates that the effect of the NN architectures on the quality of the solution depends on the structure of the trial solutions. Moreover, this demonstrates the importance of using the additional information of solutions in the network method instead of directly using the network approach.

5. The applications also show that the method can be applied to strong nonlinear cases and more accurately detect the severe nonlinearities of the physical phenomena.

6. The randomness of the selection $D_1$ in the decomposition $D=D_1+D_2$ may be a drawback. However, Theorem 1-2 proves that the solution expressions used always exist for any selection $D_1$. Hence, the method always works well at least in the accuracy range of previous similar methods. If the operator $D_1$ is selected well, then the efficiency of the method is ensured. Because the possible selections of $D_1$ are finite for a concrete IVP of ODEs, one of the efficiencies can always be selected after several times of computations. This leads to future studies of how to use more sophisticated theories to design an ANN algorithm for solving problems of DEs.

These benefits and the ideas in this study will stimulate future studies on solving DEs' problem by using the Lie theory of DEs and ANN.

\section*{Conflict of Interest}
As far as we know, there are no conflicts of interest, financial or other conflicts between the designated author and the editors, reviewers and readers of this magazine.

\section*{Author Contributions}
\textbf{Conceptualization:} Temuer Chaolu \\
\textbf{Methodology:} Ying Wen \\
\textbf{Software:} Ying Wen \\
\textbf{Writing \-original draft:} Ying Wen \\
\textbf{Writing review \& editing:} Temuer Chaolu, Ying Wen, Xiangsheng Wang

%
%
%


\begin{thebibliography}{}

\end{thebibliography}


\begin{thebibliography}{10}

\bibitem{Rivertz2013}
Rivertz HJ.
\newblock On those ordinary differential equations that are solved exactly by
  the improved Euler method.
\newblock Archivum Mathematicum. 2013;49(1):29--34.
\newblock doi:{10.5817/AM2013-1-29}.

\bibitem{Kaps1979}
Kaps P, Rentrop P.
\newblock Generalized Runge-Kutta methods of order four with stepsize control
  for stiff ordinary differential equations.
\newblock Numerische Mathematik. 1979;33(1):55--68.

\bibitem{Noye2010}
Noye B, Rankovic M.
\newblock An accurate explicit finite diference technique for solving the
  one-dimensional wave equation.
\newblock Communications in applied numerical methods. 1986;2(6):557--561.

\bibitem{Yun2019}
Yun Ys, Wen Y, Chaolu T, Rach R.
\newblock A segmented Adomian algorithm for the boundary value problem of a
  second-order partial differential equation on a plane triangle area.
\newblock Advances in Difference Equations. 2019;2019(1):1--13.

\bibitem{Olver}
Olver PJ.
\newblock Applications of Lie groups to differential equations. vol. 107.
\newblock Springer Science \& Business Media; 2000.

\bibitem{Press}
Press WH, Press WH, Flannery BP, Teukolsky SA, Vetterling WT, Flannery BP,
  et~al.
\newblock Numerical recipes in Pascal: the art of scientific computing. vol.~1.
\newblock Cambridge university press; 1989.

\bibitem{Lagaris1998}
Lagaris IE, Likas A, Fotiadis DI.
\newblock Artificial neural networks for solving ordinary and partial
  differential equations.
\newblock IEEE transactions on neural networks. 1998;9(5):987--1000.
\newblock doi:{10.1109/72.712178}.

\bibitem{Effati2010}
Effati S, Pakdaman M.
\newblock Artificial neural network approach for solving fuzzy differential
  equations.
\newblock Information Sciences. 2010;180(8):1434--1457.

\bibitem{ASH2000}
He S, Reif K, Unbehauen R.
\newblock Multilayer neural networks for solving a class of partial
  differential equations.
\newblock Neural networks. 2000;13(3):385--396.

\bibitem{Liying2007}
Li-ying X, Hui W, Zhe-zhao Z.
\newblock The algorithm of neural networks on the initial value problems in
  ordinary differential equations.
\newblock In: 2007 2nd IEEE Conference on Industrial Electronics and
  Applications. IEEE; 2007. p. 813--816.

\bibitem{Tsoulos2009}
Tsoulos IG, Gavrilis D, Glavas E.
\newblock Solving differential equations with constructed neural networks.
\newblock Neurocomputing. 2009;72(10-12):2385--2391.

\bibitem{Beidokhti2009}
Beidokhti RS, Malek A.
\newblock Solving initial-boundary value problems for systems of partial
  differential equations using neural networks and optimization techniques.
\newblock Journal of the Franklin Institute. 2009;346(9):898--913.

\bibitem{Chakraverty2014}
Chakraverty S, Mall S.
\newblock Regression-based weight generation algorithm in neural network for
  solution of initial and boundary value problems.
\newblock Neural Computing and Applications. 2014;25(3):585--594.

\bibitem{Mall2016}
Mall S, Chakraverty S.
\newblock Application of Legendre neural network for solving ordinary
  differential equations.
\newblock Applied Soft Computing. 2016;43:347--356.

\bibitem{Dockhorn2019}
Dockhorn T.
\newblock A discussion on solving partial differential equations using neural
  networks.
\newblock arXiv preprint arXiv:190407200. 2019;.

\bibitem{LiS2020}
Li S, Wang X.
\newblock Solving ordinary differential equations using an optimization
  technique based on training improved artificial neural networks.
\newblock Soft Computing. 2021;25(5):3713--3723.

\bibitem{eweinan}
Weinan E, Han J, Jentzen A.
\newblock Deep learning-based numerical methods for high-dimensional parabolic
  partial differential equations and backward stochastic differential
  equations.
\newblock Communications in Mathematics and Statistics. 2017;5(4):349--380.

\bibitem{2021Data}
L.~Wang, Z.~Yan, Data-driven peakon and periodic peakon travelling wave
  solutions of some nonlinear dispersive equations via deep learning.

\bibitem{Tongyong}
Hornik K, Stinchcombe M, White H.
\newblock Multilayer feedforward networks are universal approximators.
\newblock Neural networks. 1989;2(5):359--366.

\bibitem{aa1}
Raissi M, Perdikaris P, Karniadakis GE.
\newblock Physics informed deep learning (part i): Data-driven solutions of
  nonlinear partial differential equations.
\newblock arXiv preprint arXiv:171110561. 2017.

\bibitem{2020Solving}
Z.~Zhou, Z.~Yan, Solving forward and inverse problems of the logarithmic
  nonlinear schrodinger equation with pt-symmetric harmonic potential via deep
  learning, Physics Letters A.

\bibitem{Dis}
Dissanayake M, Phan-Thien N.
\newblock Neural-network-based approximations for solving partial differential
  equations.
\newblock communications in Numerical Methods in Engineering.
  1994;10(3):195--201.

\bibitem{Bluman}
Bluman GW, Cheviakov AF, Anco SC.
\newblock Construction of Mappings Relating Differential Equations.
\newblock In: Applications of Symmetry Methods to Partial Differential
  Equations. Springer; 2010. p. 121--186.

\bibitem{Ibragimov}
Ovsyannikov LV.
\newblock Lectures on the theory of group properties of differential equations.
\newblock World Scientific Publishing Company; 2013.

\bibitem{chao}
Chaolu T, Bluman G.
\newblock An algorithmic method for showing existence of nontrivial
  non-classical symmetries of partial differential equations without solving
  determining equations.
\newblock Journal of Mathematical Analysis and Applications.
  2014;411(1):281--296.

\bibitem{Grobner}
Gr{\"o}bner W, Knapp H.
\newblock Contributions to the method of Lie series. vol. 802.
\newblock Bibliographisches Institut Mannheim; 1967.

\bibitem{Shokin}
Shokin II, Ianenko N.
\newblock The method of differential approximation: application to gas
  dynamics.
\newblock Novosibirsk Izdatel Nauka. 1985;.

\bibitem{Doron}
Dorodnitsyn V.
\newblock Transformation groups in net spaces.
\newblock Journal of Soviet mathematics. 1991;55(1):1490--1517.

\bibitem{Sun}
Sun JQ, Ma ZQ, Hua W, Qin MZ.
\newblock New conservation schemes for the nonlinear Schr{\"o}dinger equation.
\newblock Applied mathematics and computation. 2006;177(1):446--451.

\bibitem{perturbation}
Nayfeh AH.
\newblock Introduction to perturbation techniques.
\newblock John Wiley \& Sons; 2011.

\end{thebibliography}
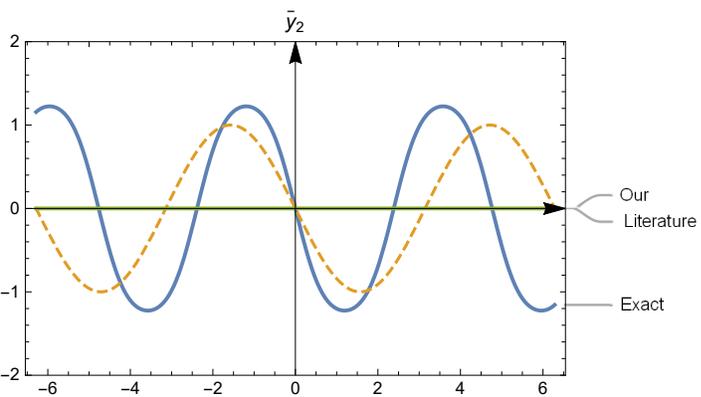
\end{document}